\newfont{\Bbb}{msbm10 scaled\magstephalf}
\documentclass[11pt]{amsart}
\usepackage{amsmath, amsthm, amscd, amsfonts}
\newtheorem{thm}{Theorem}
 
 \newtheorem{lem}[thm]{Lemma}

 \numberwithin{equation}{section}

 \def\nqp{\frac{n+1+q}{p}}
 \def\nz{(1-|z|^2)}
 \def\na{(1-|a|^2)}
 \def\nza{|1-<z,a>|}
 \def\be{\begin{equation}}
 \def\ee{\end{equation}}
 \def\bea{\begin{equation}\begin{array}{rcl}}
 \def\eea{\end{array}\end{equation}}

\begin{document}
\title[EXTENDED CES$\acute{A}$RO OPERATORS]{{EXTENDED CES$\acute{A}$RO OPERATORS  BETWEEN GENERALIZED
BESOV SPACES AND BLOCH TYPE SPACES IN THE UNIT BALL} }
\author[Z.H.Zhou and M.Zhu]{Zehua Zhou$^*$ \and  Min Zhu}
\address{\newline Department of Mathematics\newline
Tianjin University
\newline Tianjin 300072\newline P.R. China.}

\email{zehuazhou2003@yahoo.com.cn}

\address{\newline Department of Mathematics\newline
Tianjin University
\newline Tianjin 300072\newline P.R. China.}
\email{loveminjie426@tom.com}



\date{}
\thanks{\noindent $^*$ Zehua Zhou, Corresponding author. Supported in part by the National Natural Science Foundation of
China (Grand Nos.10671141, 10371091), and LiuHui Center for Applied
Mathematics, Nankai University \& Tianjin University.}

\begin{abstract}
Let $g$ be a holomorphic map of $B$, where $B$ is the unit ball of
${C}^n$. Let $0<p<+\infty, -n-1<q<+\infty$, $q>-1$ and $\alpha>0$.
This paper gives some necessary and sufficient conditions for the
Extended Ces$\acute{a}$ro Operators induced by $g$ to be bounded or
compact between generalized Besov space $B(p,q)$ and $\alpha$- Bloch
space ${\mathcal B}^\alpha.$
\end{abstract}

\maketitle

\section{Introduction}

Let $f(z)$ be a holomorphic function  on the unit disc $D$ with
$Taylor$ expansion $f(z)=\sum\limits^{\infty}_{j=0}a_jz^j$, the
classical Ces$\acute{a}$ro operator acting on $f$ is
$${\mathcal C}[f](z)=\sum\limits^{\infty}_{j=0}\left(\frac{1}{j+1}\sum\limits^j_{k=0}a_k\right)z^j.$$

In the recent years, boundedness and compactness of extended
Ces$\acute{a}$ro operator between several spaces of holomorphic
functions have been studied by many mathematicians. It is known that
the operator ${\mathcal C}$ is bounded on the usual Hardy spaces
$H^p(D)$ for $0<p< \infty$. Basic results facts on Hardy spaces can
be found in \cite{Du}. For $1\leq p<\infty$, Siskais \cite{Si1}
studied the spectrum of ${\mathcal C}$, as a by-product he obtained
that ${\mathcal C}$ is bounded on $H^p(D)$. For $p=1$, the
boundedness of ${\mathcal C}$ was given also by Sisakis \cite{Si3}
by a particularly elegant method, independent of spectrum theory, a
different proof of the result can be found in \cite{GM}. After that,
for $0<p<1$, Miao \cite{Mia} proved ${\mathcal C}$ is also bounded.
For $p=\infty$, the boundedness of ${\mathcal C}$ was given by
Danikas and Siskais in \cite{DS}. It has been also shown that the
operator ${\mathcal C}$ ia also bounded on the Bergman space (in
\cite{Si4}) as well as on the weighted Bergman spaces (in \cite{AS}
and \cite {BC}). But the operator ${\mathcal C}$ is not always
bounded, in \cite{SR}, Shi and Ren  gave a sufficient and necessary
condition for the operator ${\mathcal C}$ to be bounded on mixed
norm spaces in the unit disc.

The generalized Ces$\acute{a}$ro operators ${\mathcal C}^{\gamma}$
acting on $f$ in the unit disc were first introduced in \cite{St}
and have been subsequently studied in \cite{And} and \cite{Xi}. The
adjoint operator operator of ${\mathcal C}^{\gamma}$ was considered
in \cite{And}, \cite{Gal},\cite{Si1}, \cite{St} and \cite{Xi}. Note
that when $\gamma=0,$ ${\mathcal C}^{0}={\mathcal C}.$ Stempak
proved that ${\mathcal C}^{\gamma}$ is bounded on $H^p(D)$ for
$0<p\leq 2$. For $0<p\leq 1$, his method is similar to that of Miao;
for $p=2$, it is based on the boundedness of an appropriate sequence
transformation, and an interpolation then yields the result for
$1<p<2$. After that, Andersen \cite{And} and Xiao \cite{Xi} prove
the boundedness of ${\mathcal C}^{\gamma},$ on $H^p(D)$ for $p>2$
using different methods.

More recently,  there have been many papers focused on studying the
same problems for $n$-dimensional case, for the unit polydisc, we
refer the reader to see \cite{CS}, where they prove the boundedness
of the generalized Ces$\acute{a}$ro operator on Hardy space
$H^p(D^n)$ and the generalized Bergman space.

Let $dv$ be the $Lebesegue$ measure on the unit ball $B$ of $C^n$
normalized so that $v(B)=1$. $H(B)$ is the class of all holomorphic
functions on $B$.

A little calculation shows ${\mathcal
C}[f](z)=\frac{1}{z}\int^z_0f(t)(\log\frac{1}{1-t})'dt$. From this
point of view, if $g\in H(B)$, it is natural to consider the
extended Ces$\acute{a}$ro operator $T_g$ on $H(B)$ defined by
$$
T_g(f)(z)=\int^1_0f(tz)g(tz)\frac{dt}{t},$$ where $f\in H(B), z\in
B.$

It is easy to show that $T_g$ take $H(B)$ into itself. In general,
there is no easy way to determine when a extended Ces$\acute{a}$ro
operator is bounded or compact.

Motivated by [ShiR], Hu gave some sufficient and necessary
conditions for the extended ${\mathcal C}$ to be bounded and compact
on mixed norm spaces, Bloch space as well as Dirichlet space in the
unit ball (see [Hu1],[Hu2] and [Zha]).

For $a\in B$, let $g(z,a)= \log|\varphi_a(z)|^{-1}$ be the {\sl
Green's} function on B with logarithmic singularity at $a$, where
$\varphi_a$ is the $M\ddot{o}bius$ transformation of $B$ with
$\varphi_a(0)=a, \varphi_a(a)=0, \varphi_a=\varphi_a^{-1}$.

Let $0<p,s<+\infty, -n-1<q<+\infty$ and $q+s>-1$. We say $f\in
F(p,q,s)$ provided that $f\in H(B)$ and $$ \|
f\|_{F(p,q,s)}=|f(0)|+\{\sup_{a\in B}\int_B|\nabla
f(z)|^p(1-|z|^2)^q g^s(z,a)dv(z)\}^{\frac{1}{p}}<+\infty ,
$$
where
$$
\nabla f(z)=(\frac{\partial f(z)}{\partial z_1},\cdots,
\frac{\partial f(z)}{\partial z_n}),
$$

$F(p,q,s)$ is defined first by \cite{zhao}, we also refer the reader
to see \cite{ZhoCh}.

Let $0<p<+\infty$, $-n-1<q<+\infty$ and $q>-1$. We say $f\in B(p,q)$
provided that $f\in H(B)$ and
$$\|f\|_{(p,q)}=\{\int_B|\nabla f(z)|^p(1-|z|^2)^q dv(z)
\}^{\frac{1}{p}}<+\infty  \label{1},$$ where
$$
\nabla f(z)=(\frac{\partial f(z)}{\partial z_1},\cdots,
\frac{\partial f(z)}{\partial z_n}).
$$

It is obvious that $B(p,q)=F(p,q,0)$ if we take $s=0$. In fact,
$B(p,q)$ is also classical Besov space if we take special parameters
of $p,q$. It is not hard to show that is a $Banach$ space under the
norm $\|f\|_{B(p,q)}=|f(0)|+\|f\|_{(p,q)}$, we refer the reader to
see Zhu's book \cite{Zhu1}. From Exercises 2.2 in [Zhu1] we know
that a holomorphic function $f\in B(p,q)$ if and only if
$\int_B|Rf(z)|^p(1-|z|^2)^q <+\infty$, where $Rf(z)=<\nabla
f(z),\bar{z}>=\sum\limits^n_{j=1}z_j\frac{\partial f(z)}{\partial
z_j}.$

For $\alpha \geq 0$, $f$ is said to be in the $Bloch$ space
${\mathcal B}^\alpha$ provided that $f\in H(B)$ and $$
\|f\|_{\alpha}=\sup_{z\in B}(1-|z|^2)^\alpha|\nabla f(z)|
<+\infty.$$ As we all know, ${\mathcal B}^\alpha$ is a $Banach$
space when $\alpha\geq 1$ under the norm $\|f\|_{{\mathcal
B}^{\alpha}}=|f(0)|+\|f\|_{\alpha}$. The spaces ${\mathcal B}^1$ and
${\mathcal B}^\alpha (0<\alpha<1)$ are just the $Bloch$ space and
the $Lipschitz$ spaces $L_{1-\alpha}$ respectively. From [YaOuy] we
know that a holomorphic function $f\in{\mathcal B}^\alpha$ if and
only if $\sup_{z\in B}(1-|z|^2)^\alpha|Rf(z)|<+\infty.$ \label{3}

Furthermore, by the Norm Equivalent Theorem we have $$
\|f\|_{{\mathcal B}^\alpha}\approx |f(0)|+\sup_{z\in
B}(1-|z|^2)^\alpha|Rf(z)|, \label{4} $$ where $M \approx N$ means
the two quantities $M$ and $N$ are comparable, that is there exist
two positive constants $C_1$ and $C_2$ such that $C_1 M\leq N\leq
C_2 M$.

For $p>0, z\in B,$ denote the function $${G_p}(z)=\left\{
\begin{array}{ll}
1, & 0<p<1;\\
\log\frac{2}{1-|z|^2} ,  & p=1; \\
\Big{(}\frac{1}{1-|z|^2}\Big{)}^{\alpha-1},& p>1.
\end{array}\right.$$

In this paper, we discussed the extended Ces$\acute{a}$ro operator
between the generalized Besov space $B(p,q)$ and Bloch type space
${\mathcal B}^{\alpha}$ on the unit ball, and gave some sufficient
and necessary conditions for the operator to be bounded and compact.
The main results of the paper are the following:

\textbf{Theorem 1.}   $0<p<+\infty$, $-n-1<q<+\infty$, $q>-1,
\alpha\geq 0$, $g\in H(B)$, $T_g$ is bounded from $B(p,q)$ to
${\mathcal B}^{\alpha}$ if and only if
$$\sup_{z\in B}(1-|z|^2)^{\alpha}G_{\nqp}(z)|Rg(z)|<\infty.
$$

\textbf{Theorem 2.}   For $0<p<+\infty, -n-1<q<+\infty$, $q>-1,
\alpha\geq 0$, $g\in H(B) $, $T_g$ is compact from $B(p,q)$ to
${\mathcal B}^{\alpha}$ if and only if

 (1) If $ 0 < \nqp < 1 $, $g\in
{\mathcal B}^{\alpha};$

(2) If $\nqp \leq 1$, $\lim_{|z|\rightarrow
1^{-}}(1-|z|^2)^{\alpha}G_{\nqp}(z)|Rg(z)|= 0.$

\section{Some Lemmas}

In the following, we will use the symbol $c$ or $C$ to denote a
finite positive number which does not depend on variable $z$ and may
depend on some norms and parameters $p, q, n, \alpha, x, f$ etc, not
necessarily the same at each occurrence.

In order to prove the main result, we will give some Lemmas first.

\begin{lem} If $0<p<+\infty$, $-n-1<q<+\infty$, $q>-1$, then
$B(p,q)\subset{\mathcal B}^\nqp$ and $\exists$ $c>0$ s.t. for
$\forall f\in B(p,q)$,
$$\|f\|_{{\mathcal B}^\nqp}\leq c\|f\|_{B(p,q)}.
$$
\end{lem}

{\bf Proof.}  Suppose $f \in B(p,q).$ Fixed $0<r_0<1,$ since
$(Rf)\circ\varphi_a\in H(B),$ so $|(Rf)\circ\varphi_a|^p$ is
subharmonic in $B$. That is
\begin{eqnarray*}
&&|Rf(a)|^p =|(Rf)\circ\varphi_a(0)|^p\\
&&\leq\frac{1}{r_0^{2n}}\int_{r_0B}|(Rf)\circ\varphi_a(\omega)|^pdv(\omega)\\
&&=\frac{1}{r_0^{2n}}\int_{\varphi_a(r_0B)}|(Rf(z))|^p\frac{\na^{n+1}}{\nza^{(2n+2)}}dv(z).
\end{eqnarray*}
 From (5) in [ZhuOuy], we have
$$
\frac{1-r_0}{1+r_0}\na\leq\nz\leq\frac{1+r_0}{1-r_0}\nz
$$
as $z\in\varphi_a(r_0 B)$. Thus $$
\frac{\na^{n+1}}{\nza^{2n+2}\nz^q}\leq\frac{4^{n+1}}{\na^{n+1+q}}(\frac{1+r_0}{1-r_0})^{|q|}.
$$

 Therefore, we get
$$
\begin{array}{ll}
&|Rf(a)|^p\leq\frac{1}{r_0^{2n}}\int_{\varphi_a(r_0B)}|Rf(z)|^p\frac{\na^{n+1}}{\nza^{2n+2}}dv(z)\\
&=\frac{1}{r_0^{2n}}\int_{\varphi_a(r_0B)}|Rf(z)|^p\nz^q\frac{\na^{n+1}}{\nza^{2n+2}\nz^q}dv(z)\\
&\leq\frac{4^{n+1}r_0^{-2n}}{\na^{n+1+q}}(\frac{1+r_0}{1-r_0})^{|q|}\|f\|^p_{B(p,q)}.
\end{array}
$$
This shows that $f\in{\mathcal B}^{\frac{n+1+q}{p}}$ and
$\|f\|_{{\mathcal B}^\nqp}\leq c\|f\|_{B(p,q)}. $

\begin{lem}Let $p>0$, then there is a constant $c>0, $
for $\forall f\in{\mathcal B}^p$ and $\forall z \in B,$ the estimate
$$ |f(z)|\leq c G_p(z)\|f\|_{{\mathcal B}^p},  $$ holds,
where $G_p(z)$ is the function defined in Introduction.
\end{lem}

{\bf Proof.} This Lemma can be easily obtained by some integral
estimates. For the convenience of the reader, we will still give the
proof here.

For  $\forall f\in\beta^p(B_n),$ since
$||f||_{\beta^p}=|f(0)|+\sup\limits_{z\in B_n}(1-|z|^2)^p|\nabla
f(z)|$, we have
$$
|f(0)| \leq ||f||_{\beta^p},\quad \mbox{and}\quad |\nabla f(z)| \leq
\frac{||f||_{\beta^p}}{(1-|z|^2)^p}.
$$
but
$$
f(z)=f(0)+ \int_0^1<z,\overline{\nabla f(tz)}>dt.
$$
therefore
\begin{eqnarray*}
&&|f(z)|\leq |f(0)|+\int_0^1|z|\,|\nabla f(tz)|dt \\[6pt]
&\leq& ||f||_{\beta^p}+||f||_{\beta^p}\int_0^1
\frac{1}{(1-|tz|^2)^p}dt \leq
||f||_{\beta^p}\Big(1+\int_0^{|z|}\frac{dt}{(1-t^2)^p}\Big).
\end{eqnarray*}
when $p=1$,
$\displaystyle\int_0^{|z|}\displaystyle\frac{dt}{1-t^2}=\displaystyle\frac{1}{2}
\ln\displaystyle\frac{1+|z|}{1-|z|}\leq
\frac{1}{2}\displaystyle\ln\frac{4}{1-|z|^2}$, therefore
$$
|f(z)|\leq
\Big(1+\frac{1}{2}\ln\frac{4}{1-|z|^2}\Big)||f||_{\beta^p}.
$$

If $p\neq 1$, then
$$
\displaystyle\int_0^{|z|}\displaystyle\frac{dt}{(1-t^2)^p}=\displaystyle\int_0^{|z|}
\displaystyle\frac{dt}{(1-t)^p(1+t)^p}\leq
\displaystyle\int_0^{|z|}\displaystyle\frac{dt}{(1-t)^p}=\displaystyle\frac{1-(1-|z|)^{1-p}}{1-p},
$$
therefore when $0<p<1$, notice that
$\displaystyle\int_0^{|z|}\displaystyle\frac{dt}{(1-t^2)^p} \leq
\displaystyle\frac{1}{1-p} $ we get
$$
|f(z)|\leq \Big(1+\frac{1}{1-p}\Big)||f||_{\beta^p}.
$$
and when $p>1$
\begin{eqnarray*}
&&\int_0^{|z|}\frac{dt}{(1-t^2)^p}\leq
\frac{1-(1-|z|)^{1-p}}{1-p}\\
&&=\frac{1-(1-|z|)^{p-1}} {(p-1)(1-|z|)^{p-1}}\leq
\frac{2^{p-1}}{(p-1)(1-|z|^2)^{p-1}}
\end{eqnarray*}
so
$$
|f(z)|\leq
\Big(1+\frac{2^{p-1}}{(p-1)(1-|z|^2)^{p-1}}\Big)||f||_{\beta^p}.
$$

\begin{lem}Let $0<p<1$, $\{f_j\}$ is any bounded sequence in ${\mathcal B}^p$ and
$f_j(z)\rightarrow 0$ on any compact subset of $B$. Then
$$\lim_{j\rightarrow \infty}\sup_{z\in B}|f_j(z)|=0.$$
\end{lem}

{\bf Proof.} This lemma has been given by [Zha].

\begin{lem} There is a constant $c>0$ such that for $\forall$ $t>-1$ and $z\in
B$, $$
\int_B|\log\frac{1}{1-<z,w>}|^2\frac{(1-|w|^2)^t}{(1-<z,w>)^{n+1+t}}dv(w)\leq
C\big{(}\log\frac{1}{1-|z|^2}\big{)}^{2}.$$
\end{lem}

{\bf Proof.} This Lemma can be proved by Stirling formula and some
complex integral estimates. For the convenience of the reader, we
will still give the proof here.

Denote the right term as $I_t$ and let $2\lambda=t+n+1$. By Taylor
expansion
$$
 |\log\frac{1}{1-<z,w>}|^2=\sum\limits_{u,v=1}^{+\infty}\frac{<z,w>^u<w,z>^v}{uv}
$$
and
$$
\frac{1}{|1-<z,w>|^{2\lambda}}=\sum\limits_{k,l=0}^{+\infty}\frac{\Gamma(\lambda+k)\Gamma(\lambda+l)}{k!l!\Gamma(\lambda)^2}<z,w>^k<w,z>^l,
$$
therefore
$$
\begin{array}{ll}
I_t&=\int_B\sum\limits_{u,v=1}^{+\infty}\sum\limits_{k,l=0}^{+\infty}\frac{\Gamma(\lambda+k)\Gamma(\lambda+l)}{uvk!l!\Gamma(\lambda)^2}<z,w>^{k+u}<w,z>^{l+v}(1-|w|^2)^tdv(w)\\
   &=\sum\limits_{u=1}^{+\infty}\sum\limits_{k=0}^{+\infty}\sum\limits_{l=0}^{u+k-1}\frac{\Gamma(\lambda+k)\Gamma(\lambda+l)}{u(u+k-l)k!l!\Gamma(\lambda)^2}\int_B|<z,w>|^{2(u+k)}(1-|w|^2)^tdv(w)
\end {array}
$$
without lost of generality, let $z=|z|e_1,$ then
$$
\begin{array}{ll}
\int_B&|<z,w>|^{2(u+k)}(1-|w|^2)^tdv(w)\\
&       =\int_B (|z|w_1)^{2(u+k)}(1-|w|^2)^tdv(w)\\
                                &=2n\int_0^1\int_{\partial B}\rho^{2n-1}|z|^{2(u+k)}|\rho\xi_1|^{2(u+k)}(1-\rho^2)^td\rho d\delta_n(\xi)\\
          &=2n|z|^{2(u+k)}\int_0^1\rho^{2(u+k+n-1)+1}(1-\rho^2)^t d\rho \int_{\partial B}|\xi_1|^{2(u+k)}d\delta(\xi)\\
          &=n|z|^{2(u+k)}\frac{\Gamma(u+k+n)\Gamma(t+1)}{\Gamma(u+k+n+t+1)}\frac{(n-1)!(u+k)!}{(u+k+n-1)!}\\
          &=\frac{\Gamma(t+1)\Gamma(u+k+1)n!}{\Gamma(2\lambda+u+k)}|z|^{2(u+k)},
\end {array}
$$
so
$$
\begin{array}{ll}
I_t&=\sum\limits_{u=1}^{+\infty}\sum\limits_{k=0}^{+\infty}\sum\limits_{l=0}^{u+k-1}\frac{\Gamma(\lambda+k)\Gamma(\lambda+l)}{u(u+k-l)k!l!\Gamma(\lambda)^2}\frac{\Gamma(t+1)\Gamma(u+k+1)n!}{\Gamma(2\lambda+u+k)}|z|^{2(u+k)}\\
   &=\sum\limits_{u=1}^{+\infty}\sum\limits_{k=0}^{+\infty}\frac{n!\Gamma(t+1)\Gamma(\lambda+k)\Gamma(u+k+1)}{uk!\Gamma(\lambda)^2\Gamma(2\lambda+u+k)}\sum\limits_{l=0}^{u+k-1}\frac{\Gamma(\lambda+l)}{(u+k-l)l!}|z|^{2(u+k)}\\
   &=\sum\limits_{u=1}^{+\infty}\sum\limits_{k=1}^{+\infty}\frac{n!\Gamma(t+1)\Gamma(\lambda+k)\Gamma(u+k+1)}{uk!\Gamma(\lambda)^2\Gamma(2\lambda+u+k)}\sum\limits_{l=0}^{u+k-1}\frac{\Gamma(\lambda+l)}{(u+k-l)l!}|z|^{2(u+k)}\\
   &\ +\sum\limits_{u=1}^{+\infty}\frac{n!\Gamma(t+1)\Gamma(u+1)}{u\Gamma(\lambda)\Gamma(2\lambda+u)}\sum\limits_{l=0}^{u-1}\frac{\Gamma(\lambda+l)}{(u-l)l!}|z|^{2u}\\
   &=I_1+I_2 ,
\end{array}
$$
by Stirling formula, there is an absolute constant $C_1$
s.t.$$\frac{\Gamma(\lambda+l)}{l!}\leq C_1 l^{\lambda-1},
\frac{\Gamma(u+k+1)}{\Gamma(2\lambda +u+k)}\leq
C_1(u+k)^{1-2\lambda}, $$$$ \frac{\Gamma(u+k+1)}{\Gamma(2\lambda
+u)}\leq C_1u^{1-2\lambda}, \frac{\Gamma(\lambda+k)}{k!}\leq C_1
k^{\lambda-1}$$ for all $l,u,k\geq 1$, then
$$
I_1\leq
C_1^3\sum\limits_{u=1}^{+\infty}\sum\limits_{k=1}^{+\infty}\frac{n!\Gamma(t+1)k^{\lambda-1}(u+k)^{1-2\lambda}}{u\Gamma(\lambda)^2}\sum\limits_{l=1}^{u+k-1}\frac{l^{\lambda-1}}{(u+k-l)}|z|^{2(u+k)}
$$
and
$$
I_2 \leq
C_1^2\sum\limits_{u=1}^{+\infty}\frac{n!\Gamma(t+1)u^{1-2\lambda}}{u\Gamma(\lambda)}\sum\limits_{l=1}^{u-1}\frac{l^{\lambda-1}}{(u-l)}|z|^{2u}
.
$$
Notice that
$$
\sum\limits_{l=1}^{M-1}\frac{l^{(\lambda-1)}}{M-l}\approx
M^{\lambda-2}\log M
$$
for any $M\geq 2$, then there is constant $C$, s.t.
$$
\begin{array}{ll}
I_1\leq  &C\sum\limits_{u=1}^{+\infty}\sum\limits_{k=1}^{+\infty} \frac{n!\Gamma(t+1)k^{\lambda-1}(u+k)^{1-2\lambda}}{\Gamma(\lambda)^2u}(u+k)^{\lambda-2}\log(u+k)|z|^{2(u+k)}\\
          &=C\sum\limits_{u=1}^{+\infty}\sum\limits_{k=1}^{+\infty}\frac{n!\Gamma(t+1)}{\Gamma(\lambda)^2}\frac{k^\lambda}{(u+k)^\lambda}\frac{\log(u+k)}{u+k}\frac{1}{uk}|z|^{2(u+k)}\\
          &\leq C\sum\limits_{u=1}^{+\infty}\sum\limits_{k=1}^{+\infty}\frac{1}{uk}|z|^{2(u+k)}=C\big{(}\log\frac{1}{1-|z|^2}\big{)}^{2}
\end{array}
$$
and
$$
\begin{array}{ll}
I_2\leq  &C\sum\limits_{u=1}^{+\infty} \frac{n!\Gamma(t+1)u^{1-2\lambda}}{\Gamma(\lambda)u}u^{\lambda-2}\log u|z|^{2u}\\
          &=C\sum\limits_{u=1}^{+\infty}\frac{n!\Gamma(t+1)}{\Gamma(\lambda)}\frac{1}{u^{\lambda+1}}\frac{\log
          u}{u}|z|^{2u},
\end{array}
$$
then it is clearly that $I_2$ can be control by
$\big{(}\log\frac{1}{1-|z|^2}\big{)}^{2}$. This ends the proof of
the lemma.

\begin{lem} Let $g$ be a holomorphic self-map of $B$, $K$ is an
arbitrary compact subset of $B$. Then $T_g : B(p,q) \rightarrow
{\mathcal B}^{\alpha}$ is compact if and only if for any uniformly
bounded sequence $\{f_j\} (j\in N)$ in $B(p,q)$ which converges to
zero uniformly for $z$ on $K$ when $j\rightarrow \infty$,
$\|T_gf_j\|_{{\mathcal B}^{\alpha}}\rightarrow 0$ holds.\end{lem}

{\bf Proof.} Assume that $T_g$ is compact and suppose $\{f_j\}$ is a
sequence in $B(p,q)$ with $\sup_{j\in N}\|f_j\|_{B(p,q)}<\infty$ and
$f_j\rightarrow 0$ uniformly on compact subsets of $B$. By the
compactness of $T_g$ we have that $\{T_gf_j\}$ has a subsequence
$\{T_gf_{j_m}\}$ which converges in $\beta^{\alpha}$, say, to h. By
Lemma 2 we have that for any compact set $K\subset B$, there is a
positive constant $C_K$ independent of $f$ such that
$$
|T_gf_j(z) - h(z)|\leq C_K\|T_gf_j - h\|_{\beta^{\alpha}}
$$
for all $z\in K$. This implies that $T_gf_j(z) - h(z)\rightarrow 0$
uniformly on compact sets of $B$. Since $K$ is a compact subset of
$B$, by the hypothesis and the definition of $T_g$, $T_gf_j(z)$
converges to zero uniformly on $K$. It follows from the arbitrary of
$K$ that the limit function $h$ is equal to $0$. Since it's true for
arbitrary subsequence of $\{f_j\}$, we see that $T_gf_j\rightarrow
0$ in $\beta^{\alpha}$.

Conversely, $\{f_j\}\in K_r = B_{B(p,q)}(0,r)$, where
$B_{B(p,q)}(0,r)$ is a ball in $B(p,q)$, then by Lemma 2, $\{f_j\}$
is uniformly bounded in arbitrary compact subset $M$ of $B$. By
$Montel's$ Lemma, $\{f_j\}$ is a normal family , therefore there is
a subsequence $\{f_{j_m}\}$ which converges uniformly to $f\in H(B)$
on compact subsets of $B$. It follows that $\nabla f_{j_m}\to\nabla
f$ uniformly on compact subsets of B.

Denote $B_k=B(0,1-\frac{1}{k})\subset C^n$, then
$$\begin{array}{ll}
&\int_B|\nabla f|^p\nz^q dv(z)\\
&=\lim\limits_{k\to+\infty}\int_{B_k}\lim\limits_{m\to+\infty}|\nabla
f_{j_m}|^p\nz^qdv(z)\\
&\leq\lim\limits_{k\to+\infty}\lim\limits_{m\to+\infty}\int_{B_k}|\nabla
f_{j_m}|^p\nz^q dv(z).
\end{array}
$$
But $\{f_{j_m}\}\subset B_{B(p,q)}(0,r)$, then $$\int_{B_k}|\nabla
f_{j_m}|^p\nz^q dv(z)<r^p,$$therefore
$$
\int_B|\nabla f|^p\nz^q dv(z)\leq r^p.
$$
So $\|f\|_{B(p,q)}\leq r,$ and $f\in B(p,q).$ Hence the sequence
$\{f_{j_m}-f\}$ is such that $\|f_{j_m}-f\|\leq 2r<\infty$ and
converges to $0$ on compact
subsets of B, by the hypothesis of this lemma, we have that
$$T_gf_{j_m} \to T_gf$$ in ${\mathcal B}^{\alpha}$. Thus the set $T_g(K_r)$ is relatively compact,
finishing the proof.

\begin{lem} Let $g\in H(B)$, then $$
R[T_gf](z)=f(z)Rg(z)$$ for any $f\in H(B)$ and $z\in B$.
\end{lem}

{\bf Proof.}  Suppose the holomorphic function $fRg$ has the
$Taylor$ expansion $$(f R g)(z)=\sum_{|\alpha|\geq
1}a_{\alpha}z^{\alpha}.$$ Then we have
\begin{eqnarray*}
&&R(T_gf)(z)=R\int^1_0 f(tz) R(tz) \frac{dt}{t}
=R\int^1_0\sum_{|\alpha|\geq 1}a_{\alpha}(tz)^{\alpha}\frac{dt}{t}\\
&&=R [\sum_{|\alpha|\geq 1}\frac{a_{\alpha}z^{\alpha}}{|\alpha|}]
=\sum_{|\alpha|\geq 1}a_{\alpha}z^{\alpha}=(f R g)(z).
\end{eqnarray*}

\section{The Proof Of Theorem 1}

Suppose $\sup_{z\in B}(1-|z|^2)^{\alpha}G_{\nqp}(z)|Rg(z)|<\infty $.
$\forall f\in H(B) $ then by Lemma 1, Lemma 2 and Lemma 6, we have
\begin{eqnarray*}
&&(1-|z|^2)^{\alpha}| R [T_gf](z) |\\
&&=(1-|z|^2)^{\alpha} |f(z)| |Rg(z)|\\
&&\leq c (1-|z|^2)^{\alpha}G_{\nqp}(z) |Rg(z)|\\
&&\leq c \|f\|_{B(p,q)}(1-|z|^2)^{\alpha}G_{\nqp}(z) |Rg(z)|\\
&&\leq c \|f\|_{B(p,q)}.
\end{eqnarray*}

Therefore, $T_g$ is  bounded .

On the other hand, suppose  $T_g $  is  bounded,
with$$\|T_gf\|_{{\mathcal B}^{\alpha}} \leq c\|f\|_{B(p,q)}.$$

  (1) If  $0 < \nqp < 1$ , it's very easy to show that the function
$f(z)=1$ are in $B(p,q)$, therefore $T_gf$ must be in
${\mathcal B}^{\alpha}$,\\
namely
\begin{eqnarray*}
&&\sup_{z\in B}(1-|z|^2)^{\alpha}|RT_gf(z)|\\
&&=\sup_{z\in B}(1-|z|^2)^{\alpha}|Rg(z)| < \infty.
\end{eqnarray*}

 (2) If  $\nqp > 1$ , we need to prove that  $\sup_{z\in
  B}(1-|z|^2)^{\alpha}(\frac{1}{1-|z|^2})^{\nqp-1}|R g(z)|< \infty $.

For $w \in B$, take the test function
$$f_w(z)=\frac{1-|w|^2}{(1-<z,w>)^{\nqp}}.$$

It is easy to see that
$$\int_B (1-|z|^2)^q|\nabla f_w(z)|^p dv(z)
\leq c(1-|w|^2)^p \int_B \frac{(1-|z|^2)^q}{|1-<z,w>|^{n+1+q+p}}
dv(z)\leq c.$$

The last inequality follws from [Zhu1], so $f_w \in B(p,q)$ for any
$w\in B$. With the boundedness of $T_g$, we get
\begin{eqnarray*}
&&(1-|z|^2)^{\alpha}(\frac{1}{1-|z|^2})^{\nqp-1}|Rg(z)|\\
&&=(1-|z|^2)^{\alpha}|f_z(z)||Rg(z)|\\
&&=(1-|z|^2)^{\alpha}|R(T_gf_z)(z)|\\
&& \leq\|T_gf_z\|_{{\mathcal B}^{\alpha}} \leq c\|T_g\| < \infty.
\end{eqnarray*}

  (3) If  $\nqp = 1$, namely $p=n + 1 + q$, we need to prove
  $$\sup_{z\in B}
(1-|z|^2)^{\alpha}\log\frac{2}{1-|z|^2}|Rg(z)| < \infty.$$

For $ w \in B$, take the test function
$$f_w(z)=(\log\frac{1}{1-|w|^2})^{-\frac{2}{p}}(\log\frac{1}{1-<z,w>})^{1+\frac{2}{p}}.$$

It is easy to show that $f_w\in B(p,q)$ from Lemma 4. The same
discussion as the case (2) gives the needed result, and we omit it
here. So, the proof of Theorem 1 is completed.

\section{The Proof Of Theorem 2 }

 $\{f_j\}$ is an uniformly bounded sequence in $B(p,q)$ which
converges to zero uniformly on any compact subset of $B$ when
$j\rightarrow \infty$.

(1) If $T_g$ is compact, we have got that $g\in {\mathcal
B}^{\alpha}$.

On the other hand, from Lemma 1, we know that $\|f_j\|_{{\mathcal
B}^{\nqp}} \leq c\|f_j\|_{B(p,q)}$, thus $\{f_j\}$ is
 unformly  bounded in ${\mathcal B}^{\nqp}$. Then by the hypothesis and  Lemma 3, we get
 that $$\lim_{j\rightarrow \infty}\sup_{z\in B}|f_j(z)|=0.$$

 Therefore
$$
 \|T_gf_j\|_{{\mathcal B}^{\alpha}}
\leq c\sup_{z\in B}(1-|z|^2)^{\alpha}|f_j(z)Rg(z)|
 \leq c\|g\|_{{\mathcal B}^{\alpha}}\sup_{z\in B}|f_j(z)|.
$$

 Then when $j\rightarrow \infty$, $\|T_gf_j\|_{{\mathcal B}^{\alpha}} \rightarrow
 0$. So $T_g$ is compact from Lemma 5.

(2) If $\lim_{|z|\rightarrow
1^{-}}(1-|z|^2)^{\alpha}G_{\nqp}(z)|Rg(z)|= 0$, then $\forall
 \varepsilon >0$, $\exists  r \in  (0,1)$, such that $$
(1-|z|^2)^{\alpha}G_{\nqp}(z)|Rg(z)| < \varepsilon ,    r<|z|<1.
$$

Then
\begin{eqnarray*}
&&\|T_gf_j\|_{{\mathcal B}^{\alpha}} \leq  c\sup_{|z|\leq
r}(1-|z|^2)^{\alpha}|f_j(z)Rg(z)| + c\sup_{r<|z|<1}(1-|z|^2)^{\alpha}|f_j(z)Rg(z)|\\
&&\leq c\sup_{|z| \leq r}(1-|z|^2)^{\alpha}|Rg(z)| |f_j(z)| + c
\sup_{r<|z|<1}(1-|z|^2)^{\alpha}G_{\nqp}(z)|Rg(z)| \|f_j\|_{B(p,q)}\\
&&\leq c\sup_{|z| \leq r}(1-|z|^2)^{\alpha} |Rg(z)| |f_j(z)| + c \varepsilon  \|f_j\|_{B(p,q)}\\
&&\leq c \varepsilon ,
\end{eqnarray*}

if $j$ is sufficiently large. This means $ \|T_gf_j\|_{{\mathcal
B}^{\alpha}} \rightarrow  0$ as $j$ tends to $\infty$.

On the other hand, if $\nqp = 1$, it is sufficient to prove
$$\lim_{|z|\rightarrow
1}(1-|z|^2)^{\alpha}|Rg(z)|\log\frac{1}{1-|z|^2} = 0.$$

Suppose that $\lim_{|z|\rightarrow
1}(1-|z|^2)^{\alpha}|Rg(z)|\log\frac{1}{1-|z|^2}\neq 0$, then there
exists  $\varepsilon_0 >0 $, $\{z^j\}\in B$, such that
$$(1-|z^j|^2)^{\alpha}|Rg(z^j)|\log\frac{1}{1-|z^j|^2} \geq
\varepsilon_0.
$$

Let
$$f_j(z)=(\log\frac{1}{1-|z^j|^2})^{-\frac{2}{p}}(\log\frac{1}{1-<z,z^j>})^{1+\frac{2}{p}}.$$

We have shown that $f_j \in B(p,q)$  with $\|f_j\|_{B(p,q)} \leq c $
, and it is obvious that $f_j\rightarrow 0 $ uniformly on any
compact subset of $B$ as $j \rightarrow \infty$. While
\begin{eqnarray*}
&& \|T_gf_j\|_{{\mathcal B}^{\alpha}}\\
&& \geq (1-|z^j|^2)^{\alpha}|f_j(z^j)||R g(z^j)|\\
&& =\{(1-|z^j|^2)^{\alpha}|Rg(z^j)|\log\frac{1}{1-|z^j|^2}\}
 |f_j(z^j)| (\log\frac{1}{1-|z^j|^2})^{-1}\\
&& \geq \varepsilon_0 |f_j(z^j)| (\log\frac{1}{1-|z^j|^2})^{-1}\\
&& = \varepsilon_0,
\end{eqnarray*}

then $\|T_gf_j\|_{{\mathcal B}^{\alpha}}$ doesn't tend to $0$ when
$j \rightarrow \infty$. It's a contraction. So
$$\lim_{|z|\rightarrow
1}(1-|z|^2)^{\alpha}|Rg(z)|\log\frac{1}{1-|z|^2} = 0.
$$

Meanwhile, as  $ \lim_{|z|\rightarrow 1}\log\frac{1}{1-|z|^2} =
\infty $, it is easy to see that $\lim_{|z| \rightarrow
1}(1-|z|^2)^{\alpha}|Rg(z)| = 0$.

Therefore, we have $$\lim_{|z|\rightarrow
1}(1-|z|^2)^{\alpha}|Rg(z)|\log\frac{2}{1-|z|^2} = 0 .$$

If $\nqp > 1$, just let $$
f_j(z)=\frac{1-|z^j|^2}{(1-<z,z^j>)^{\nqp}} ,$$

and use the same method as the situation of  $ \nqp = 1 $, we can
also prove that the theorem holds. So, the proof of Theorem 2 is
completed.

\end{document}